\documentclass{amsart}
\usepackage{amssymb}
\usepackage{amscd}
\usepackage{enumerate}
\usepackage{graphicx}
\usepackage{url}

\newtheorem{theorem}{Theorem}[section]
\newtheorem{lemma}[theorem]{Lemma}

\newtheorem{proposition}[theorem]{Proposition}

\newtheorem{question}[theorem]{Question}

\theoremstyle{definition}
\newtheorem{definition}[theorem]{Definition}

\newtheorem{remark}[theorem]{Remark}

\title{Gromov hyperbolicity and a variation of the Gordian complex}

\author{Kazuhiro Ichihara}
\address{School of Mathematics Education,
Nara University of Education,
Takabatake-cho, Nara, 630-8528, Japan}
\email{ichihara@nara-edu.ac.jp}
\urladdr{http://mailsrv.nara-edu.ac.jp/~ichihara/index.html}
\thanks{The first author is partially supported by
Grant-in-Aid for Young Scientists (B), No. 20740039,
Ministry of Education, Culture, Sports, Science and Technology, Japan.}

\author{In Dae Jong}
\address{Graduate School of Science, Osaka City University, Osaka 558-8585, Japan}
\email{jong@sci.osaka-cu.ac.jp}
\urladdr{http://www.ex.media.osaka-cu.ac.jp/~d07sa009/index.html}

\date{\today}



\begin{document}

\maketitle

\begin{abstract}
We introduce new simplicial complexes by using various invariants and local moves for knots, which give generalizations of the Gordian complex defined by Hirasawa and Uchida. 
In particular, we focus on the simplicial complex defined by using the Alexander-Conway polynomial and the Delta-move, and show that the simplicial complex is Gromov hyperbolic and quasi-isometric to the real line. 
\end{abstract}


\section{Introduction}\label{sect:intro}

\begin{figure}[h]
\begin{center}
\includegraphics[width=0.12\textwidth]{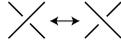}
\caption{The crossing change.}\label{fig:cc}
\end{center}
\end{figure}

A {\it knot} is an ambient isotopy class of a simple closed curve smoothly embedded in the $3$--sphere. 
Let $\mathcal{K}$ be the set of all knots. 
Let $\lambda$ be a local move on knots, that is a local operation deforming a knot (see~\cite[Section 2]{OhyamaYamada2008} for the precise definition of a local move). 
The {\it $\lambda$-Gordian distance} $d^\lambda(K,K')$ between knots $K$ and $K'$ is defined to be the minimal number of the local moves $\lambda$ needed to deform $K$ into $K'$. 
If such a minimum does not exist, then we set $d^{\lambda}(K,K') = \infty$. 
Let $\mathrm{x}$ denote the {\it crossing change} which is a local move as shown in Figure~\ref{fig:cc}. 
In the case where $\lambda = \mathrm{x}$, $d^\mathrm{x}$ is called the {\it Gordian distance}. 

Using the Gordian distance, Hirasawa and Uchida~\cite{HirasawaUchida2002} defined the {\it Gordian complex} denoted by $\mathcal{G}^\mathrm{x}$. 
More generally, the {\it $\lambda$-Gordian complex} $\mathcal{G}^\lambda$ introduced by Nakanishi and Ohyama~\cite[Section 1]{NakanishiOhyama2006} is the simplicial complex defined by the following; 
\begin{itemize}
\item the set of vertices (i.e.~$0$--simplices) of $\mathcal{G}^\lambda$ is $\mathcal{K}$, and 
\item $n+1$ vertices $K_0, \dots, K_n$ span an $n$--simplex if and only if $d^\lambda(K_i, K_j) = 1$ holds for each $i \ne j \in \{\, 0, \dots, n\,\}$. 
\end{itemize}
We call the $1$--skelton of $\mathcal{G}^\mathrm{x}$ (resp.~$\mathcal{G}^\lambda$) the {\it Gordian graph} (resp.~the {\it $\lambda$-Gordian graph}), and denote it by $G^\mathrm{x}$ (resp.~$G^\lambda$). 
Assuming that every edge has length $1$, each connected component of $G^\lambda$ is regarded as a metric space which turns to a geodesic space (see Section~\ref{sec:pre}). 
Then one of the problems we are interested in is to reveal properties on such spaces. 
In particular, we are interested in global properties, and the Gromov hyperbolicity~\cite{Gromov1987} (for a brief review, see Section~\ref{sec:pre}) is an important one. 
There are several studies on simplicial complexes arising in geometry and topology, in particular, the curve complex introduced by Harvey~\cite{Harvey1981} is widely studied (see also \cite{Hamentadt2007} for a survery). 
Masur and Minsky proved that the curve complex is Gromov hyperbolic~\cite{MasurMinsky1999}. 
On the other hand, there is no known fact on the Gromov hyperbolicity of $G^\lambda$ except for the following. 

\begin{proposition}[{\cite[Theorem C]{GambaudoGhys2005}}]
The Gordian graph $G^\mathrm{x}$ is not Gromov hyperbolic.
\end{proposition}

\smallskip

Here we will introduce new simplicial complexes and graphs by using knot invariants and local moves, which give generalizations of the $\lambda$-Gordian complex and the $\lambda$-Gordian graph. 
Let $\iota$ be a knot invariant, that is, a function on $\mathcal{K}$ such that $\iota(K)$ and $\iota(K')$ coincide if $K$ is equivalent to $K'$. 
We denote by $K \sim_\iota K'$ if $\iota(K) = \iota(K')$ holds. 
Clearly the binary relation $\sim_\iota$ provides an equivalence relation on $\mathcal{K}$. 
Let $[K]_\iota$ denote the equivalence class of $K$, and $\mathcal{K}_\iota = \{\, [K]_\iota \, | \, K \in \mathcal{K} \,\}$. 

\begin{definition}\label{def:vlG}
Let $\iota$ be a knot invariant, and let $\lambda$ be a local move on knots. 
The {\it $(\iota,\lambda)$-Gordian complex} $\mathcal{G}^\lambda_\iota$ is defined by the following; 
\begin{itemize}
\item The set of vertices of $\mathcal{G}^\lambda_\iota$ is $\mathcal{K}_\iota$, and  
\item $n+1$ vertices $[K_0]_\iota, \dots, [K_n]_\iota$ span an $n$--simplex if and only if for each $i \ne j \in \{\, 0, \dots, n\,\}$, there exists a pair of knots $K_{i,j} \in [K_i]_\iota$ and $K_{j,i} \in [K_j]_\iota$ such that $d^\lambda (K_{i,j}, K_{j,i}) = 1$. 
\end{itemize}
\end{definition}

We call the $1$--skelton of $\mathcal{G}^\lambda_\iota$, denoted by $G^\lambda_\iota$, the {\it $(\iota,\lambda)$-Gordian graph}. 
Assuming that every edge has length $1$, we regard $G^\lambda_\iota$ as a metric space. 
Then we denote by $d_\iota^\lambda$ the metric on $G^\lambda_\iota$, and call it the {\it $(\iota,\lambda)$-Gordian distance}. 

\smallskip

\begin{figure}
\begin{center}
\includegraphics[height=0.035\textheight]{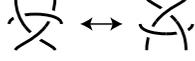}
\caption{The Delta-move.}\label{fig:DeltaMove}
\end{center}
\end{figure}

\begin{figure}
\begin{center}
\includegraphics[height=0.035\textheight]{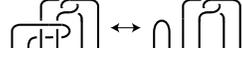}
\caption{A $C_2$-move.}\label{fig:C2}
\end{center}
\end{figure}

Let $\nabla_K$ be the Conway polynomial~\cite{Conway1970} of a knot $K$, which is a polynomial in $z^2$ with integer coefficients. 
The Conway polynomial is also called the Alexander-Conway polynomial since it is regarded as a normalized Alexander polynomial~\cite{Alexander1928}. 
We refer the reader to \cite{KawauchiBook} for basic terminologies of knot theory. 
The {\it Delta-move}, denoted by the symbol $\Delta$, is a local move on knots as shown in Figure~\ref{fig:DeltaMove}, which was introduced by Matveev~\cite{Matveev1987}, and Murakami and Nakanishi~\cite{MurakamiNakanishi1989} independently. 
It is known that the Delta-move is equivalent to a $C_2$-move (see Figure~\ref{fig:C2}), which is one of $C_n$-moves introduced by Goussarov~\cite{Goussarov1995} and Habiro~\cite{Habiro2000} independently. 

Using the Conway polynomial and the Delta move, the $(\nabla,\Delta)$-Gordian graph $G^\Delta_\nabla$ are defined. 
In this paper, we show the following. 

\begin{theorem}\label{thm:main}
The $(\nabla,\Delta)$-Gordian graph $G^\Delta_\nabla$ is $2$-hyperbolic. 
Further it is quasi-isometric to the real line $\mathbb{R}$. 
\end{theorem}

\begin{remark}\label{rem:see}
In Section~\ref{sec:Q} we will see that $G^\Delta_\nabla$ and $\mathcal{G}^\Delta_\nabla$ coincide (see Proposition~\ref{pro:1dim}). 
Thus, $\mathcal{G}^\Delta_\nabla$ is also $2$-hyperbolic and quasi-isometric to $\mathbb{R}$. 
\end{remark}

This paper is constructed as follows: 
In Section~\ref{sec:pre}, we give a brief review on a Gromov hyperbolic space to study $G^\Delta_\nabla$. 
In Section~\ref{sec:distance}, we study the $(\nabla,\Delta)$-Gordian distance. 
In Section~\ref{sec:proof}, we prove Theorem~\ref{thm:main}. 
In Section~\ref{sec:Q}, we give observations on the complexes $\mathcal{G}^\mathrm{x}_\nabla$ and $\mathcal{G}^\Delta_\nabla$. 
We also give some remarks and questions related to our study. 


\subsection*{Acknowledgments}
The authors would like to thank Masayuki Asaoka, Yoshifumi Matsuda, Yasutaka Nakanishi, Ken'ichi Ohshika, Yoshiyuki Ohyama, and Harumi Yamada for their various suggestions and comments.

\section{Preliminaries}\label{sec:pre}

In this section, we give a brief review on a Gromov hyperbolic space. 
For details, see \cite{BridsonHaefliger1999} or \cite{Gromov1987}. 
Let $X$ be a geodesic space, that is, a metric space such that the distance between any two points is equal to the length of a geodesic segment joining them. 
We denote by $s(x,y)$ a geodesic segment joining two points $x$ and $y$. 
A {\it geodesic triangle} $T$ in $X$ is a triple of points $x, y, z \in X$ together with three geodesic segments $s(x,y)$, $s(y,z)$, and $s(z,x)$ called the {\it sides} of $T$. 
For $\delta \ge 0$, a geodesic triangle is said to be {\it $\delta$-slim} if each side of a triangle belongs to the $\delta$-neighborhood of the union of the other two sides. 
We say that $X$ is {\it $\delta$-hyperbolic} (or {\it Gromov hyperbolic}) if there exists a constant $\delta \ge 0$ such that any geodesic triangle in $X$ is $\delta$-slim. 
Clearly if a geodesic space is $\delta$-hyperbolic for a particular $\delta$, then it is also $\delta'$-hyperbolic for all $\delta' \ge \delta$. 

Let $\Gamma$ be a connected graph. 
We denote by $\overline{vv'}$ an edge adjacent to vertices $v$ and $v'$. 
In general, assuming that each edge has length $1$, the graph $\Gamma$ is regarded as a metric space, and it turns to a geodesic space. 

As mentioned in Section~\ref{sect:intro}, we regard $G^\Delta_\nabla$ as a geodesic space with the metric $d^\Delta_\nabla$ induced by the above setting. 
Note that $G^\Delta_\nabla$ is a connected graph since the Delta-move is an unknotting operation~\cite{MurakamiNakanishi1989}. 
Recall that the vertex set of $G^\Delta_\nabla$ is $\mathcal{K}_\nabla$. 
After now, unless otherwise specified, we use the symbol $[K]$ to refer the vertex $[K]_\nabla$ for brevity. 
Note that for two vertices $[K]$ and $[K']$ with $d^\Delta_\nabla([K],[K']) = p$, a geodesic segment $s([K],[K'])$ is of the form $\overline{v_0v_1} \cup \overline{v_1v_2} \cup \cdots \cup \overline{v_{p-1}v_p}$ for some vertices $v_0, v_1, \ldots,v_p \in \mathcal{K}_\nabla$ with $v_0 = [K]$ and $v_p = [K']$.

Let $X$ and $X'$ be metric spaces with metric functions $d$ and $d'$ respectively. 
A map $f : X \rightarrow X'$ is a {\it quasi-isometry} if there exist constants $A, E \ge 0$, $B,C,D >0$  such that $A d(x,y) - B \le d'(f(x),f(y)) \le C d(x,y) + D$ holds for any $x,y \in X$, and for any $x' \in X'$ there exists $x \in X$ such that $d'(x',f(x)) \le E$. 
Then we say that $X$ is {\it quasi-isometric} to $X'$. 
It is known that the Gromov hyperbolicity is an invariant for geodesic spaces under quasi-isometries.

\section{$(\nabla,\Delta)$-Gordian distance}\label{sec:distance}

In this section, we study the $(\nabla,\Delta)$-Gordian distance. 
Let $a_n(K)$ be the $n$-th coefficient of $\nabla_K$. 
The following lemma is a fundamental relationship between the Conway polynomials and the Delta-move. 

\begin{lemma}[\cite{Okada1990}]\label{lem:Okada}
For $K, K' \in \mathcal{K}$ with $d^\Delta(K,K') = 1$, we have $a_2(K) - a_2(K') = \pm1.$
Furthermore, for any $[K], [K'] \in \mathcal{K}_\nabla$, we have
$$d_\nabla^\Delta([K],[K']) \ge |a_2(K) - a_2(K')|\, ,$$ 
and 
$$d_\nabla^\Delta([K],[K']) \equiv |a_2(K) - a_2(K')| \mod 2 \,.$$ 
\end{lemma}

Note that $a_2(K_1) = a_2(K_2)$ holds for any $K_1 ,K_2 \in [K]$. 
The following lemma gives the formula to detect the $(\nabla, \Delta)$-Gordian distance between any pair of vertices. 

\begin{lemma}\label{lem:DeltaDistance}
For $[K] \ne [K'] \in \mathcal{K}_\nabla$, we have the following. 
\begin{align*}
d_\nabla^\Delta([K],[K']) = 
\begin{cases}
|a_2(K) - a_2(K')| & a_2(K) \ne a_2(K')\\
2 & a_2(K) = a_2(K')
\end{cases}
\end{align*}
\end{lemma}
\begin{proof}
It is known that every vertex $[K]$ contains an unknotting number one knot~\cite{Kondo1979}, \cite{Sakai1977}. 
Let $J \in [K]$ be an unknotting number one knot. 
Then for any sequence of integers $(m_4,\dots,m_{2l})$, there exists a knot $J'$ satisfying $d^\Delta(J,J') = 1$, $a_2(J) - a_2(J') = 1$, and $a_{2j}(J) - a_{2j}(J') = m_{2j}$ for $j=2,\dots, l$~\cite[Lemma A]{NakanishiOhyama2006}. 
\begin{itemize}
\item 
If $a_2 (K) \ne a_2(K')$, then by the above argument, we have $d^\Delta_\nabla([K],[K']) \le |a_2(K) - a_2(K')|$. 
On the other hand, by Lemma~\ref{lem:Okada}, we have $d^\Delta_\nabla([K],[K']) \ge |a_2(K) - a_2(K')|$. 
Hence we have $d^\Delta_\nabla([K],[K']) = |a_2(K) - a_2(K')|$. 
\item
If $a_2 (K) = a_2(K')$, then by Lemma~\ref{lem:Okada} and the assumption $[K] \ne [K']$, $d^\Delta_\nabla ([K],[K'])$ is a non-zero even integer, namely $d^\Delta_\nabla ([K],[K']) \ge 2$. 
On the other hand, by the above argument, we have $d^\Delta_\nabla ([K],[K']) \le 2$. 
Hence we have $d^\Delta_\nabla([K],[K']) = 2$. 
\end{itemize}
Now we complete the proof of Lemma~\ref{lem:DeltaDistance}. 
\end{proof}

\section{Proof of Theorem~\ref{thm:main}}\label{sec:proof}

For $\varepsilon \ge 0$, let $N(p,\varepsilon)$ be the $\varepsilon$-neighborhood of a point $p \in G^\Delta_\nabla$, and $N(P,\varepsilon)$ the $\varepsilon$-neighborhood of a subset $P \subset G^\Delta_\nabla$, that is, $N(p,\varepsilon) = \{\, q \in G^\Delta_\nabla \,|\, d^\Delta_\nabla(p,q) \le \varepsilon \,\}$ and $N(P,\varepsilon) = \bigcup_{p \in P} N(p,\varepsilon)$. 
Let $V_n = \{\, [K] \in \mathcal{K}_\nabla \,|\, a_2(K) = n \,\}$. 
Then we have the following. 

\begin{lemma}\label{lem:3nbd}
For any $[K] \in \mathcal{K}_\nabla$ with $a_2(K) = n$, we have $$N([K],2) \supset N(V_{n},1) \, .$$ 
\end{lemma}
\begin{proof}
Recall that a {\it vertex-induced subgraph} is a subset of the vertices together with any edges whose endpoints are both in this subset. 
Then we see that $N(V_{n},1)$ is the vertex-induced subgraph which is induced by the subset $V_{n-1} \cup V_{n} \cup V_{n+1}$ of vertices. 
By Lemma~\ref{lem:DeltaDistance}, we have $d^\Delta_\nabla ([K],v_{n}) = 2$ for any $v_n \ne [K] \in V_n$, and $d^\Delta_\nabla([K], v_{n \pm1}) = 1$ for any $v_{n \pm1} \in V_{n \pm 1}$. 
This completes the proof of Lemma~\ref{lem:3nbd}. 
\end{proof}

Now we start the proof of Theorem~\ref{thm:main}. 

\begin{figure}[h]
\begin{center}
\includegraphics[width=0.5\textwidth]{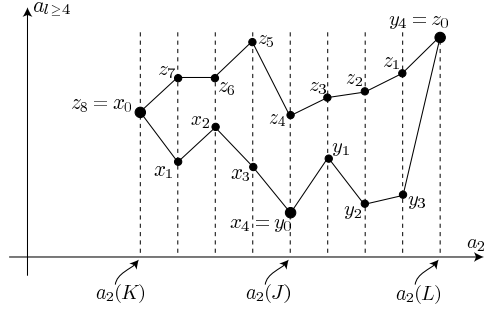}
\caption{The case where $p=4$ and $q=4$.}\label{fig:position}
\end{center}
\end{figure}

\begin{proof}[Proof of Theorem~\ref{thm:main}]
Let $T$ be a geodesic triangle in $G^\Delta_\nabla$ with sides $s(x,y)$, $s(y,z)$, and $s(z,x)$. 
We only show the case where $x$, $y$, and $z$ are in $\mathcal{K}_\nabla$ since the other cases (i.e.~the cases where some of $x$, $y$, and $z$ are not contained in $\mathcal{K}_\nabla$) are proved in a similar way. 
Let $x = [K]$, $y = [J]$, and $z = [L]$. 
Without loss of generality, we may assume that $a_2(K) \le a_2(J) \le a_2(L)$. 
Let $k = a_2(J) - a_2(K)$ and $k' = a_2(L) - a_2(J)$. 
Let 
\begin{align*}
s(x,y) &= \overline{x_0x_1} \cup \overline{x_1x_2} \cup \cdots \cup \overline{x_{p-1}x_{p}},\\ 
s(y,z) &= \overline{y_0y_1} \cup \overline{y_1y_2} \cup \cdots \cup \overline{y_{q-1}y_{q}},\\ 
s(z,x) &= \overline{z_0z_1} \cup \overline{z_1z_2} \cup \cdots \cup \overline{z_{r-1}z_{r}}, 
\end{align*}
where $x_0, \dots, x_p, y_0,\dots,y_q, z_0,\dots,z_r$ are in $\mathcal{K}_\nabla$ with $x_0 = x = z_r$, $y_0 = y = x_p $, and $z_0 = z = y_q$. 
We show that $T$ is $\delta$-slim for a positive integer $\delta \le 2$, that is, $G^\Delta_\nabla$ is $2$-hyperbolic.

\begin{description}

\item[Case 1. {\boldmath $k \ge 1$}, {\boldmath $k' \ge 1$}] 
By Lemma~\ref{lem:DeltaDistance}, we have $p = k $, $q = k'$, and $r = k+k'$. 
Figure~\ref{fig:position} is an example of a geodesic triangle for $p=q=4$. 
(In Figure~\ref{fig:position}, we plot vertices with respect to the coefficients of the Conway polynomial.) 
First we show that $N(s(x,y)\cup s(y,z),2) \supset s(z,x)$. 
By Lemma~\ref{lem:3nbd}, we have 
$$N(y_j,2) \supset \overline{z_{q-j+1}z_{q-j}}\, , \, \overline{z_{q-j}z_{q-j-1}}$$ 
for each $j = 1, \cdots , q-1$. 
Thus, we have 
\begin{align*}
N(s(y,z),2) &\supset N(y_1 \cup \cdots \cup y_{q-1},2)\\
&\supset \overline{z_{q}z_{q-1}} \cup \cdots \cup \overline{z_{1}z_{0}}. 
\end{align*}
Similarly, we have 
\begin{align*}
N(s(x,y),2) &\supset N(x_1 \cup \cdots \cup x_{p-1},2)\\
&\supset \overline{z_{r}z_{r-1}} \cup \cdots \cup \overline{z_{q+1}z_{q}}. 
\end{align*}
Therefore we have $N(s(x,y)\cup s(y,z),2) \supset s(z,x)$. 
Remaining two conditions $N(s(y,z)\cup s(z,x),2) \supset s(x,y)$ and $N(s(z,x)\cup s(x,y),2) \supset s(y,z)$ are shown by the similar argument. 
Therefore the geodesic triangle $T$ is $2$-slim.

\item[Case 2. {\boldmath $k=0$}, {\boldmath $k' \ge 1$}] 
By Lemma~\ref{lem:DeltaDistance}, we have $p=2$, $q = r= k'$. 
Then by Lemma~\ref{lem:3nbd}, we see that $T$ is $2$-slim.

\item[Case 3. {\boldmath $k \ge 1$}, {\boldmath $k'=0$}] 
This case is proved by the same argument applied in Case 2.

\item[Case 4. {\boldmath $k=0$}, {\boldmath $k'=0$}] 
By Lemma~\ref{lem:DeltaDistance}, we have $p = q = r = 2$. 
Then by Lemma~\ref{lem:3nbd}, we see that $T$ is $1$-slim. 

\end{description}
Therefore $G^\Delta_\nabla$ is $2$-hyperbolic. 

Next we show that $G^\Delta_\nabla$ is quasi-isometric to the real line $\mathbb{R}$. 
Let $f : G^\Delta_\nabla \rightarrow \mathbb{R}$ be a map defined by the following; 
$f(\overline{v_{n-1}v_{n}}) = [n-1,n]$ for $v_{n-1} \in V_{n-1}$ and $v_n \in V_n$. 
Here $[n-1,n]$ denotes the closed interval bounded by $n-1$ and $n$. 
Then, by Lemma~\ref{lem:DeltaDistance}, the map $f$ is a quasi-isometry. 
\end{proof}

\section{Known facts and questions}\label{sec:Q}

In this section, we introduce some facts and remarks related to the diameter and the dimension of complexes, and propose questions for further studies. 

First we focus on the $(\nabla,\Delta)$-Gordian complex $\mathcal{G}^\Delta_\nabla$. 
By Lemma~\ref{lem:DeltaDistance}, there exists no $2$--simplex in $\mathcal{G}^\Delta_\nabla$ (see {\cite[Proposition 2.3]{Ohyama2006}}). 
This implies the following. 
\begin{proposition}\label{pro:1dim}
The dimension of $\mathcal{G}^\Delta_\nabla$ is one. 
Therefore $\mathcal{G}^\Delta_\nabla$ and $G^\Delta_\nabla$ coincide. 
\end{proposition}

\smallskip

Next we focus on the $(\nabla,\mathrm{x})$-Gordian complex $\mathcal{G}^\mathrm{x}_\nabla$ and the $(\nabla,\mathrm{x})$-Gordian graph $G^\mathrm{x}_\nabla$. 
As mentioned in the proof of Lemma~\ref{lem:DeltaDistance}, any Alexander-Conway polynomial of a knot is realized by an unknotting number one knot~\cite{Kondo1979}, \cite{Sakai1977}. 
(There are several studies on the realization problem of the Alexander-Conway polynomial. 
See \cite{Fujii1996}, \cite{Jong2009, Jong2009a}, \cite{Levine1965}, \cite{Nakamura2009}, \cite{RolfsenBook}, \cite{Seifert1935}.) 
Thus, the diameter of $\mathcal{G}^\mathrm{x}_\nabla$ and that of $G^\mathrm{x}_\nabla$ are less than or equal to two. 
Since a geodesic space with a finite diameter $r$ is $r$-hyperbolic, and it is quasi-isometric to a point, we have the following. 
\begin{proposition}\label{pro:NablaX}
The $(\nabla,\mathrm{x})$-Gordian graph $G^\mathrm{x}_\nabla$ is $2$-hyperbolic, and it is quasi-isometric to a point.  
\end{proposition}

\begin{figure}[h]
\begin{center}
\includegraphics[width=0.16\textwidth]{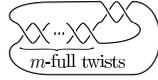}
\caption{A twist knot $K_m$.}\label{fig:TwistKnot}
\end{center}
\end{figure}

\begin{remark}\label{rem:diam-infty}
For any $n \in \mathbb{N}$, there exists vertices $[K]$ and $[K']$ such that $d^\Delta_\nabla([K],[K']) = n$. 
Actually, we have $d^\Delta_\nabla([K_0],[K_n]) = n$, where $K_0$ and $K_n$ are twist knots depicted in Figure~\ref{fig:TwistKnot}. 
Thus, the diameter of $G^\Delta_\nabla$ is infinite. 
\end{remark}

Recently, Kawauchi~\cite{Kawauchi2009} showed that $d^\mathrm{x}_\nabla ([K_1], [K_{-1}]) = 2$ by using duality theorems on the infinite cyclic covering space of a knot exterior. 
Here $K_1$ denotes the trefoil knot and $K_{-1}$ denotes the figure-eight knot as shown in Figure~\ref{fig:TwistKnot}. 
Thus, the diameter of $G^\mathrm{x}_\nabla$ is just two. 
Kawauchi also showed that the dimension of the $(\nabla,\mathrm{x})$-Gordian complex $\mathcal{G}^\mathrm{x}_\nabla$ is infinite~\cite{Kawauchi2009}, that is, an $n$--simplex is contained in $\mathcal{G}^\mathrm{x}_\nabla$ for any $n \in \mathbb{N}$.

\begin{remark}\label{rem:dimension}
The dimension of Gordian complex $\mathcal{G}^\mathrm{x}$ is infinite, which was shown by Hirasawa and Uchida~\cite{HirasawaUchida2002}. 
For a $C_n$-move with $n \ge 3$, the dimension of the $C_n$-Gordian complex $\mathcal{G}^{C_n}$ is also infinite, which was shown by Ohyama~\cite{Ohyama2006}. 
Here we note that the $C_1$-move is the crossing change and a $C_2$-move is equivalent to the Delta-move.  
\end{remark}

Finally we propose some questions. 
For several local moves and knot invariants, it is interesting to consider the Gromov hyperbolicity of $G^\lambda_\iota$. 
As mentioned in Remark~\ref{rem:dimension}, the crossing change is equivalent to the $C_1$-move and the Delta-move is equivalent to a $C_2$-move. 
Then the following question is natural to ask. 

\begin{question}\label{que:CnComplex}
For $n \ge 3$, is each connected component of the $(\nabla, C_n)$-Gordian graph $G^{C_n}_\nabla$ Gromov hyperbolic? 
\end{question}

Note that $G^{C_n}_\nabla$ consists of infinitely many connected components by results of Goussarov~\cite{Goussarov1995} and Habiro~\cite{Habiro2000}. 

It is also interesting to study the $\lambda$-Gordian graph. 
In particular, $G^\mathrm{x}_\nabla$ and $G^\Delta_\nabla$ are Gromov hyperbolic, but $G^\mathrm{x}$ is not Gromov hyperbolic. 

\begin{question}\label{que:DeltaComplex}
Is the $\Delta$-Gordian graph $G^\Delta$ Gromov hyperbolic? 
\end{question}

\providecommand{\bysame}{\leavevmode\hbox to3em{\hrulefill}\thinspace}
\providecommand{\MR}{\relax\ifhmode\unskip\space\fi MR }
\providecommand{\MRhref}[2]{%
  \href{http://www.ams.org/mathscinet-getitem?mr=#1}{#2}
}
\providecommand{\href}[2]{#2}


\begin{thebibliography}{10}

\bibitem{Alexander1928}
J.~W. Alexander, \emph{{Topological invariants of knots and links}}, Trans.
  Amer. Math. Soc. \textbf{30} (1928), 275--306.

\bibitem{BridsonHaefliger1999}
M.~Bridson and A.~Haefliger, \emph{Metric spaces of non-positive curvature},
  Grundlehren der Mathematischen Wissenschaften [Fundamental Principles of
  Mathematical Sciences], vol. 319, Springer-Verlag, Berlin, 1999.

\bibitem{Conway1970}
J.~H. Conway, \emph{An enumeration of knots and links, and some of their
  algebraic properties}, Computational Problems in Abstract Algebra (Proc.
  Conf., Oxford, 1967) (1970), 329--358.

\bibitem{Fujii1996}
H.~Fujii, \emph{{Geometric indices and the Alexander polynomial of a knot}},
  Proc. Amer. Math. Soc. \textbf{124} (1996), no.~9, 2923--2933.

\bibitem{GambaudoGhys2005}
J.-M. Gambaudo and \'E. Ghys, \emph{Braids and signatures}, Bull. Soc. Math.
  France \textbf{133} (2005), no.~4, 541--579.

\bibitem{Goussarov1995}
M.~N. Goussarov, \emph{Knotted graphs and a geometrical technique of
  $n$-equivalences}, POMI Sankt Petersburg preprint, circa (1995), (in
  Russian).

\bibitem{Gromov1987}
M.~Gromov, \emph{Hyperbolic groups}, Essays in group theory, Math. Sci. Res.
  Inst. Publ, vol.~8, pp.~75--263, Springer, New York, 1987.

\bibitem{Habiro2000}
K.~Habiro, \emph{Claspers and finite type invariants of links}, Geom. Topol.
  \textbf{4} (2000), 1--83.

\bibitem{Hamentadt2007}
U.~Hament\"adt, \emph{{Geometry of the complex of curves and of Teichm\"uller
  space}}, Handbook of Teichm\"uller Theory: I (IRMA Lectures in Mathematics \&
  Theoretical Physics), vol.~11, ch.~10, pp.~447--469, European Mathematical
  Society, 2007.

\bibitem{Harvey1981}
W.~J. Harvey, \emph{Boundary structure of the modular group}, Riemann Surfaces
  and Related Topics: Proceedings of the 1978 Stony Brook Conference (I.~Kra
  and B.~Maskit, eds.), Ann. Math. Stud., vol.~97, Princeton, 1981.

\bibitem{HirasawaUchida2002}
M.~Hirasawa and Y.~Uchida, \emph{{The Gordian complex of knots}}, J. Knot
  Theory Ramifications \textbf{11} (2002), no.~3, 363--368.

\bibitem{Jong2009}
I.~D. Jong, \emph{Alexander polynomials of alternating knots of genus two},
  Osaka J. Math. \textbf{46} (2009), no.~2, 353--371.

\bibitem{Jong2009a}
\bysame, \emph{{Alexander polynomials of alternating knots of genus two II}},
  to appear in J. Knot Theory Ramifications (2009).

\bibitem{KawauchiBook}
A.~Kawauchi, \emph{{A Survey of knot theory}}, Birckh\"{a}user-Verlag, Basel,
  1996.

\bibitem{Kawauchi2009}
A.~Kawauchi, \emph{{On the Alexander polynomials of knots with Gordian distance
  one}}, preprint (2009).

\bibitem{Kondo1979}
H.~Kondo, \emph{{Knots of unknotting number $1$ and their Alexander
  polynomials}}, Osaka J. Math. \textbf{16} (1979), no.~2, 551--559.

\bibitem{Levine1965}
J.~Levine, \emph{A characterization of knot polynomials}, Topology \textbf{4}
  (1965), 135--141.

\bibitem{MasurMinsky1999}
H.~A. Masur and Y.~N. Minsky, \emph{{Geometry of the complex of curve. I.
  Hyperbolicity}}, Invent. Math. \textbf{138} (1999), 103--149, 1.

\bibitem{Matveev1987}
S.~G. Matveev, \emph{Generalized surgeries of three-dimensional manifolds and
  representations of homology spheres}, Mat. Zametki \textbf{42} (1987), no.~2,
  268--278.

\bibitem{MurakamiNakanishi1989}
H.~Murakami and Y.~Nakanishi, \emph{On a certain move generating
  link-homology}, Math. Ann. \textbf{284} (1989), no.~1, 75--89.

\bibitem{Nakamura2009}
T.~Nakamura, \emph{{Braidzel surfaces for fibered knots with given Alexander
  polynomials}}, to appear in Kobe J. Math. (2009).

\bibitem{NakanishiOhyama2006}
Y.~Nakanishi and Y.~Ohyama, \emph{{Local moves and Gordian complexes}}, J. Knot
  Theory Ramifications \textbf{15} (2006), no.~9, 1215--1224.

\bibitem{Ohyama2006}
Y.~Ohyama, \emph{{The $C\sb k$-Gordian complex of knots}}, J. Knot Theory
  Ramifications \textbf{15} (2006), 73--80.

\bibitem{OhyamaYamada2008}
Y.~Ohyama and H.~Yamada, \emph{{A $C_n$-move for a knot and the coefficients of
  the Conway polynomial}}, J. Knot Theory Ramifications \textbf{17} (2008),
  no.~7, 771--785.

\bibitem{Okada1990}
M.~Okada, \emph{{Delta-unknotting operation and the second coefficient of the
  Conway polynomial}}, J. Math. Soc. Japan \textbf{42} (1990), no.~4, 713--717.

\bibitem{RolfsenBook}
D.~Rolfsen, \emph{{Knots and links}}, AMS Chelsea publishing, 2003.

\bibitem{Sakai1977}
T.~Sakai, \emph{{A remark on the Alexander polynomials of knots}}, Math. Sem.
  Notes Kobe Univ. \textbf{5} (1977), no.~3, 451--456.

\bibitem{Seifert1935}
H.~Seifert, \emph{{\"Uber das Geschlecht von Knoten}}, Math. Ann. \textbf{110}
  (1935), 571--592.

\end{thebibliography}

\end{document}